\font\capit=cmcsc10
\font\addressit=cmcsc8
\font\eightrm=cmr8
\font\sixrm=cmr6

\def\idfy{\mathbin{\hbox{$\widetilde\rightarrow$}}}
\def\textindent#1{\indent\llap{#1\enspace}\ignorespaces}
\def\itemitem{\par\indent \hangindent2\parindent \textindent}

\def\f2{{\bf F}_2}
\def\td{{\rm td}\ }
\def\res{{\rm res}}

\def\kum{{\rm Kum}}

\def\bl{{\rm Bl}}
\def\secant{{\rm Sec}}
\def\sectil{{\widetilde \secant}}

\def\sing{{\rm Sing}\,}

\def\ctil{{\widetilde C}}

\def\eeetil{{\widetilde \varepsilon}}

\def\dual{^{\vee}}

\def\dbar{\overline D}
\def\delbar{\overline \Delta}

\def\khat{\widehat K}

\def\peta{P_{\eta}}

\def\nm{{\rm Nm}}
\def\pic{{\rm Pic}}

\def\ff{{\cal F}}

\def\nn{{\cal N}}

\def\bb{{\cal B}}
\def\ll{{\cal L}}
\def\gg{{\cal G}}

\def\ss{{\cal S}}
\def\ll{{\cal L}}
\def\pp{{\cal P}}

\def\oo{{\cal O}}

\def\fftil{\widetilde{\cal F}}
\def\su{{\cal SU}_C}
\def\sux{{\cal SU}_X}
\def\iic{{\cal I}_C}
\def\eee{\varepsilon}

\def\gr{{\rm Gr}}
\def\grs#1#2{{\rm Grass}_{#1}(#2)}
\def\ses#1#2#3{0\rightarrow {#1}\rightarrow {#2}\rightarrow {#3}\rightarrow 0}

\def\z{{\bf Z}}

\def\c{{\bf C}}
\def\bp{{\bf P}}

\def\im{{\rm im}}

\def\ker{{\rm ker \ }}
\def\rk{{\rm rank \ }}
\def\map#1{\ \smash{\mathop{\longrightarrow}\limits^{#1}}\ }

\def\pf{{\it Proof.}\ }

\def\<{\langle}
\def\>{\rangle}

\newcommand{\qed}{{\unskip\nobreak\hfill\hbox{ $\Box$}\medskip\par}}

\documentstyle[epsf,leqno,bezier,12pt]{article}

\makeatletter                                            
\renewcommand{\@begintheorem}[2]{                        
\sl \trivlist \item [\hskip \labelsep {\bf #2\ \ #1.}]   
				}                        
\makeatother                                             
\makeatletter
\def\section{\@startsection {section}{1}{\z@}{-3.5ex plus -1ex minus
 -.2ex}{1.5ex plus .2ex}{\large\bf}}

\def\subsection{\@startsection{subsection}{2}{\z@}{-3.25ex plus -1ex
minus
 -.2ex}{1.5ex plus .2ex}{\normalsize\it}}

\newcommand{\numberequationsassubsubsections}

\newtheorem{prop}{Proposition}[section]
\newtheorem{lemm}[prop]{Lemma}
\newtheorem{theo}[prop]{Theorem}
\newtheorem{cor}[prop]{Corollary}
\newtheorem{conj}[prop]{Conjecture}
\newtheorem{rem}[prop]{\it Remark}
\newtheorem{rems}[prop]{\it Remarks}
\newtheorem{ex}[prop]{Example}

\begin{document}

\title{Heisenberg invariant quartics and $\su(2)$ for a curve of genus four}
\author{William Oxbury and Christian Pauly}
\date{}

\maketitle
\bigskip

The projective moduli variety $\su(2)$ of semistable rank 2 vector bundles with trivial determinant on a smooth projective curve $C$ comes with a natural morphism $\phi$ to the linear series $|2\Theta|$ where $\Theta$ is the theta divisor on the Jacobian of $C$. Well-known results of
Narasimhan and Ramanan say that $\phi$ is an isomorphism to $\bp^3$ if $C$ has genus 2 \cite{NR1}, and when $C$ is nonhyperelliptic of genus 3 it is an isomorphism to a special Heisenberg-invariant quartic $Q_C\subset\bp^7$ \cite{NR2}. The present paper is an attempt to extend these results to higher genus.

In the nonhyperelliptic genus 3 case
the so-called {\it Coble quartic} $Q_C \subset |2\Theta| = \bp^7$ is characterised by either of two properties:

\medskip
\itemitem{(i)}
$Q_C$ is the unique Heisenberg-invariant quartic containing the Kummer variety, i.e. the image of $\kum: J_C \rightarrow |2\Theta|$, $x\mapsto \Theta_x +\Theta_{-x}$, in its singular locus; and
  
\itemitem{(ii)} 
$Q_C$ is precisely the set of $2\Theta$-divisors containing some translate of the curve $W_1 \subset J_C^1$.

\medskip
\noindent
We shall examine, for a curve of genus 4, the analogue of each of these properties, and our first main result, analogous to (i), is:

\begin{theo}
\label{theorem1}
If $C$ is a curve of genus 4 without vanishing theta-nulls then there exists a unique (irreducible) Heisenberg-invariant quartic $Q_C\subset |2\Theta| = \bp^{15}$ containing $\phi(\su(2))$ in its singular locus.
\end{theo}

We prove this in sections \ref{invariantqs} and \ref{cubicnorm} (see corollary \ref{specialq}). The main work involved is first to 
show cubic normality
for $\phi(\su(2))$ (theorem \ref{cubicnormality}). 
We then
use the Verlinde formula to deduce that its ideal contains exactly sixteen independent cubics; by symmetry considerations these cubics are the partial derivatives, with respect to the homogeneous coordinates, of a single quartic $Q_C$. The argument here is identical to that of Coble \cite{C} for the genus 3 case. 

We conjecture that $\phi(\su(2)) = {\rm Sing}\ Q_C$, or equivalently that the ideal of $\phi(\su(2))$ is generated by cubics. We cannot prove this, but in the rest of the paper we examine the relationship of this problem with property (ii) above. For any curve one may construct a sequence of irreducible, Heisenberg-invariant subvarieties, for $1\leq d\leq g-1$,
$$
G_{d} = \{D\ |\ \hbox{$x+W_{g-d} \subset {\rm supp}\ D$ for some $x\in J_C^{d-1}$}\}
\subset |2\Theta|.
$$
In particular, $G_1$ is the Kummer variety, while
$G_{g-1}$ is a hypersurface
containing $\phi(\su(2))$ and
which coincides with the Coble quartic in the case $g=3$ (and with the Kummer quartic surface in the case $g=2$). When $g=4$, however, $G_3$ turns out to be distinct from $Q_C$---quite contrary to our original expectation.

We see this by restricting to the eigen-$\bp^7$s of the action on $|2\Theta| = \bp^{15}$ of the group $J_C[2]$ of 2-torsion points. 
For any nonzero element $\eta \in J_C[2]$ we have an associated double cover $\pi:\ctil \rightarrow C$ with respect to which $\ker \nm = \peta \cup \peta^-$, where $(\peta,\Xi)$ is the principally polarised Prym variety. 
The fixed-point set of the $\eta$-action on $|2\Theta|$ is a pair of $\bp^7$s either of which can be naturally identified with $|2\Xi|$; this fixed-point set therefore contains the Kummer image of $\peta \cup \peta^-$, and this is precisely the intersection with $\phi(\su(2))$. 

Beauville and Debarre \cite{BD} have shown that 
a $|2\Xi|$-embedded Prym-Kummer variety admits a 4-parameter family of quadrisecant planes analogous to the trisecant lines of a Jacobian Kummer. We prove:

\begin{theo}
\label{theorem2}
Let $C$ be a curve of genus 4 without vanishing theta-nulls; and for any nonzero $\eta \in J_C[2]$ identify $|2\Xi|\hookrightarrow |2\Theta|$ as the component of the fixed-point set of $\eta$ containing the Kummer image of $\peta$. Then:
\begin{enumerate}
\item
$Q_C\subset \bp^{15}$ restricts on $|2\Xi|$ to the Coble quartic of $\kum(\peta)$.
\item
$G_3\subset \bp^{15}$ restricts on $|2\Xi|$ to the hypersurface ruled by the quadrisecant planes of $\kum(\peta)$; and this is distinct from the Coble quartic. 
\end{enumerate}
\end{theo}

Part 1 is proved in section \ref{cubicnorm}. Note that $\peta$ is necessarily the Jacobian $J_X$ of some curve $X$ of genus 3, which can be constructed explicitly (given a choice of trigonal pencil on $C$) via the Recillas correspondence (see section \ref{pryms}). In section \ref{cobleqs} we obtain necessary and sufficient conditions for a secant line of $\kum(J_X)$ to lie on $\phi(\sux(2))$; and in section \ref{quadrisecants} this is used to show that a generic quadrisecant plane of the family does not lie on $\phi(\sux(2))$---hence the final remark in part 2 of the theorem. 

We introduce the filtration of $|2\Theta|$ by the subvarieties
$G_{d}$
in section \ref{abeljacobi}.
These varieties are ruled, as $x\in J_C^{d-1}$ varies, by the subseries $N(x) \subset |2\Theta|$ of divisors containing $x+W_{g-d}$. Moreover, for any theta characteristic $\kappa$ on the curve the rulings of $G_d$ and $G_{g+1-d}$ are polar in the sense that 
$$
N(x) \perp N(\kappa x^{-1})
$$
with respect to the bilinear form on $|2\Theta|$ induced by $\kappa$ (symmetric or skew-symmetric accordingly as $\kappa$ is even or odd).
Our results here are somewhat incomplete when $g>4$, and depend on the vanishing, which we are unable to prove, of certain cohomology groups on the symmetric products $S^d C$. This problem and the associated computations are discussed in the appendix.

For genus 4, however, we are able to prove what we need, and we arrive at a configuration:
$$
\kum(J) \subset G_2 \subset \phi(\su(2)) \subset G_3 \subset \bp^{15}
$$
where $G_2$ is a divisor in $\phi(\su(2))$ ruled by 4-planes, and $G_3$ is ruled by their polar 10-planes (with respect to any theta characteristic).
$G_2$ contains the trisecants of the Jacobian Kummer variety (this is true for any genus, incidentally, and is proved in \cite{OPP}); while the ruling of $G_3$ cuts out in $|2\Xi|$, for each $\eta \in J_C[2]$, precisely the Beauville--Debarre quadrisecant planes. Again this is true for any genus, and is proved in section \ref{prymG3}.

\medskip
{\it Acknowledgments.}
The authors are grateful to Arnaud Beauville and Bert van Geemen for some helpful discussions; and to Miles Reid and Warwick MRC for their hospitality during the early part of 1996, when much of this work was carried out. The second author thanks the University of Durham for its hospitality during 1996, and European networks Europroj and AGE for partial financial support.

\vfill\eject

\section{Preliminaries}
\label{prelims}

Given a smooth projective curve $C$ of genus $g\geq 2$, let $J= J^0(C)$ be its Jacobian and $\Theta = W_{g-1} \subset J^{g-1}(C)$ its canonical theta divisor. Let $\vartheta(C) \subset J^{g-1}(C)$ be the set of its theta characteristics. Then the line bundle $\ll = \oo(2\Theta_{\kappa})\in \pic(J)$ is independent of $\kappa \in \vartheta(C)$ (where $\Theta_{\kappa} = \Theta - \kappa$); and there is a canonical duality (up to scalar)
$$
H^0(J, \ll) = H^0(J^{g-1}, 2\Theta)\dual.
$$

Let $\su(2)$ and $\su(2,K)$ be the projective moduli varieties of semistable rank 2 vector bundles on $C$ with determinant $\oo_C$ and $K=K_C$ respectively;
and let $\ll$, $\ll_K$ denote the respective ample generators of their Picard groups. Then there are canonical identifications:

\begin{equation}
\label{2theta}
|\ll|_J = |2\Theta|\dual = |\ll|_{\su(2)}
\qquad
{\rm and}
\qquad
|\ll|_J\dual = |2\Theta| = |\ll_K|_{\su(2,K)}.
\end{equation}
It should not cause any great confusion that we denote by $\ll$ both the line bundle on $\su(2)$ and its pull-back via the semistable boundary $J\rightarrow \su(2)$, $L\mapsto L\oplus L^{-1}$.

These spaces give us certain maps which are all identified by (\ref{2theta}):

\begin{equation}
\begin{array}{rcl}
\phi: \su(2) &\rightarrow& |2\Theta| \\
    E &\mapsto & D_E = \{ L\in J^{g-1}\ | \ h^0(C,L\otimes E) >0\ \};\\
&&\\
 \lambda_{|\ll|} : \su(2) &\rightarrow& |\ll|\dual;\\
&&\\
\psi: \su(2) &\rightarrow& |\ll_K|\\
E &\mapsto & \{F\in \su(2,K)\ |\ h^0(C,E\otimes F) >0\ \}_{\rm red}.\\
\end{array}
\end{equation}
Note that in the definition of $\psi$, the bundle $E\otimes F$ carries a $K$-valued orthogonal structure, so by the Atiyah-Mumford lemma the condition 
$H^0(C,E\otimes F) \not=0$ determines naturally a divisor with multiplicity 2. Note also that $\phi(E) = D_E$ is the restriction of $\psi(E)$ to the semistable boundary $J^{g-1} \rightarrow \su(2,K)$, $L\mapsto L\oplus KL^{-1}$.

There are likewise (naturally identified) maps:

\begin{equation}
\begin{array}{rcl}
 \lambda_{|\ll_K|} : \su(2,K) &\rightarrow& |\ll_K|\dual;\\
&&\\
\phi_K: \su(2,K) &\rightarrow& |\ll|_J\\
F &\mapsto & \{L\in J\ |\ h^0(C,L\otimes F) >0\ \}.\\
\end{array}
\end{equation}

Note that any choice of theta characteristic $\kappa \in \vartheta (C)$ sets up a commutative diagram

\begin{equation}
\label{thetachoice}
\begin{array}{rcccl}
J&\rightarrow&\su(2)&\map{\phi}&|2\Theta| = |\ll|\dual\\
&&&&\\
\downarrow&&\otimes \kappa \downarrow&\searrow&\downarrow\\
&&&&\\
J^{g-1}&\rightarrow&\su(2,K)&\map{\phi_K}&|\ll|,\\
\end{array}
\end{equation}
in which all the vertical arrows are isomorphisms: the last $|2\Theta | \idfy |\ll|_J$ is just translation $D\mapsto D-\kappa$; let us briefly recall how the isomorphism $|\ll| \cong |\ll|\dual$ (or $|2\Theta |\cong|2\Theta |\dual$) is induced by the theta characteristic.

We shall write $V = H^0(J,\ll)$, and let $\theta_{\kappa} \in V$ denote a section cutting out the divisor $\Theta_{\kappa}$. We consider the map 
$$
\begin{array}{rcl}
m: J\times J & \rightarrow & J\times J\\
(u,v) &\mapsto & (u+v,u-v)\\
\end{array}
$$
and denote by $\xi_{\kappa} \in \bigotimes^2 V$ the element corresponding to $m^*(pr_1^*\theta_{\kappa} \otimes pr_2^*\theta_{\kappa})$ under the K\"unneth 
isomorphism
$$
\textstyle
\bigotimes^2 V \cong H^0(J\times J, pr_1^*\ll \otimes pr_2^*\ll).
$$
Then each $\xi_{\kappa}$ is a nondegenerate pairing (and we obtain bases $\{ \xi_{\kappa}\}_{\kappa \in \vartheta^+(C)}$ of $S^2 V$ and $\{ \xi_{\kappa}\}_{\kappa \in \vartheta^-(C)}$ of $\bigwedge^2 V$); this gives the isomorphism $V \cong V\dual$ on the right in diagram (\ref{thetachoice}). 
Restricted to the Kummer variety this pairing can be written

\begin{equation}
\label{riemannxi}
\xi_{\kappa} (u,v) = \theta_{\kappa}(u+v)\theta_{\kappa}(u-v)
\qquad
\hbox{for $u,v \in J$.}
\end{equation}

Finally, note that the whole diagram (\ref{thetachoice}) is acted on by the  
subgroup $J[2]$ of 2-torsion points in the Jacobian, acting on vector bundles by tensor product and on $|2\Theta|$ by translation. The horizontal maps are all equivariant and the south and south-east arrows are permuted by the action of $J[2]$ on $\vartheta(C)$.

\section{Prym varieties}
\label{pryms}

For each nonzero half-period $\eta \in J[2] \backslash
\{\oo\}$ we have an associated unramified double cover
$$\pi : \ctil \rightarrow C.$$ 
We shall denote by $\sigma$ the involution of $\ctil$ given by sheet-interchange over~$C$; and by abuse of notation it will denote also the induced involution of $\pic(\ctil)$.
The kernel of the norm map on divisors has two isomorphic connected components:
$$
\ker \nm = P_{\eta} \cup P_{\eta}^-,
$$
where $P_{\eta} = (1 - \sigma)J^0(\ctil)$ and $P_{\eta}^- = (1 - \sigma)J^1(\ctil)$. 
We shall refer to the composite
$$
S^d \ctil \rightarrow J^d(\ctil) \rightarrow P_{\eta} \cup P_{\eta}^-
$$
as the Abel-Prym map, and for $D\in S^d \ctil$
write $[D] = D-\sigma D \in P_{\eta}$ for $d$ even, $\in P_{\eta}^-$ for $d$ odd.

If we choose $\zeta \in J$ such that $\zeta^2 = \eta$ then we can define a map
$$
\peta \cup \peta^- \rightarrow \su(2)  \map{\phi} |2\Theta|
$$
by $x \mapsto \zeta \otimes \pi_* x \in \su(2)$. The image is independent of the choice of $\zeta$ and is precisely the fixed-point set of the involution $\otimes \eta$ of $\su(2)$. Moreover, the linear span of this image can be naturally identified with the linear series $|2\Xi|$ where $\Xi$ is the canonical theta divisor on the dual abelian variety, and represents the principal polarisation on $\peta$ induced from that on $J(\ctil)$. Thus there is a commutative diagram
\begin{equation}
\label{prymkum}
\begin{array}{ccc}
\peta & \rightarrow & \su(2) \\
&&\\
\downarrow &&\downarrow\\
&&\\
|2\Xi| &\subset & |2\Theta|.\\
\end{array}
\end{equation}

Finally, we recall that by a construction of S. Recillas \cite{Rec}, if $C$ is a trigonal curve then $\peta = \peta(C)$ is isomorphic as a ppav to a tetragonal Jacobian:

\begin{prop}
There is generically a bijection between the following two sets of data, under which there is a canonical isomorphism of ppavs $(J(X),\Theta) \cong (P_{\eta}(C), \Xi)$:
\begin{enumerate}
\item
$(X,g^1_4)$ where $X$ is a smooth curve of genus $g$, and $g^1_4$ is a tetragonal pencil;
\item
$(C,\eta,g^1_3)$ where $C$ is a smooth curve of genus $g+1$, $\eta \in J[2]$ is a nonzero 2-torsion point and $g^1_3$ is a trigonal pencil.
\end{enumerate}
\end{prop}

For the details of this construction we refer to \cite{D}---note in particular that the word `generically' means that on each side one must restrict to pencils having smooth ramification behaviour, in order to obtain a smooth curve on the other side. Accordingly, the correspondence can be compactified, but this will not concern us.
 
Here we need only to note that by definition the double cover defined by $\eta$ is the relative symmetric square $\ctil = S^2_{\bp^1}X$ with respect to the 4 to 1 cover $X\rightarrow \bp^1$ determined by $g^1_4$; and then $C = \ctil /\sigma$ where $\sigma $ is the obvious involution on $S^2_{\bp^1}X$. We then have a commutative diagram

\begin{equation}
\label{recdiagram}
\begin{array}{ccc}
\ctil & \hookrightarrow & S^2 X\\
&&\\
\downarrow && \downarrow \hbox{Abel-Jacobi} \hidewidth\\
&&\\
P_{\eta}^- & \idfy & J^2(X).\\
\end{array}
\end{equation}

We shall be particularly interested in the Recillas correspondence for the case $g=3$.

\section{Heisenberg invariant quartics}
\label{invariantqs}

Let $V= H^0(J,\ll) = H^0(\su(2), \ll)$; so from section \ref{prelims}
we have maps:
$$
J\rightarrow \su(2) \map{\phi}  \bp V\dual = |2\Theta|;
$$
and these identify with the maps given by the complete linear series $|\ll|$ on $J$ and $\su(2)$ respectively.

Let us fix a theta structure for the line bundle $\ll$; this allows us to 
view the vector space $V$ as an irreducible representation of the Heisenberg group
$$
H_g = \c^* \times \f2^g \times {\rm Hom}(\f2^g, \c^*),
$$$$
(s,a,\chi)(t,b,\gamma) = (st\gamma(a), a+b,\chi\gamma).
$$
This in turn gives us a canonical basis $\{X_{\sigma}\}_{\sigma \in \f2^g}$ for $V$ (i.e. homogeneous coordinates on $|2\Theta| = \bp^{2^g -1}$) such that:
$$
(s,a, \chi) : X_{b} \mapsto s\chi(a + b) X_{a + b},
\qquad
(s,a, \chi) \in H_g.
$$
In particular our theta structure fixes an isomorphism 
$$
J[2] \cong \f2^g \times {\rm Hom}(\f2^g, \c^*) = \khat \times K 
$$ 
where we denote by $K$ and $\khat$ the maximal level subgroups
$$
\begin{array}{rcl}
K = K_g &=& \{(1,0,\chi)| \chi \in {\rm Hom}(\f2^g, \c^*)\},\\
\khat = \khat_g &=& \{(1,a,1)| a \in \f2^g\}.\\
\end{array}
$$
The Heisenberg group acts also on the spaces of higher degree forms $S^nV$. This action is related to differentiation of polynomials by the following lemma, which is easy to check.

\begin{lemm}
\label{diffn}
For $x=(s,a,\chi)\in H_g$ and $P\in S^n V$ we have:
$$
{\partial \over \partial X_{b}}(x\cdot P)=
\chi(b)\
x\cdot {\partial P\over \partial X_{a+b}}.
$$
\end{lemm}

We shall be concerned, in particular,
with the spaces $S^3 V$ and $S^4 V$
of cubics and
quartics on $|2\Theta|$ respectively,
and in particular with the subspace $S^4_0 V \subset S^4 V$ of invariant quartics.
This subspace has a basis consisting of:
$$
\begin{array}{rclcl}
Q_0 &=& \sum_{\sigma \in \khat} X_{\sigma}^4,&&\\
Q_{\lambda} &=& \sum_{\sigma \in \khat} X_{\sigma}^2 X_{\sigma+\lambda}^2&&\lambda\in \bp(\khat),\\ 
Q_{\Lambda} &=& \sum_{\sigma \in \khat} X_{\sigma} X_{\sigma+\lambda} X_{\sigma+\mu} X_{\sigma+\nu}&& \Lambda = \{\lambda,\mu,\nu\}\in \grs{2}{\khat},\\
\end{array}
$$
and in particular has dimension 
\begin{equation}
\label{dimS40}
\dim S^4_0 V = {1\over 3}(2^g +1)(2^{g-1}+1).
\end{equation}

\begin{prop} 
\label{partials}
{\rm (\cite{vG}, proposition 2)}
With respect to the homogeneous coordinates $\{X_{\sigma}\}_{\sigma \in \khat}$ we have:
\begin{enumerate}
\item 
$\displaystyle
{\partial \over \partial X_0} : S^4_0 V \idfy (S^3 V )^{K}$ is an isomorphism;
\item
For $\sigma \in \khat$ and $Q \in S^4_0 V$ we have $\displaystyle{\partial Q\over \partial X_{\sigma}}= \sigma \cdot {\partial Q\over \partial X_0}$.
\end{enumerate}
\end{prop}

The second part follows at once from lemma \ref{diffn}; the proposition
says that an invariant quartic is determined by any of its 
partial derivatives with respect to the homogeneous coordinates; and these
are all $K$-invariant and are permuted by the action of $\khat$.

\begin{prop}
\label{coble}
Suppose that $M\subset |2\Theta|$ is a $J[2]$-invariant subvariety and that
the restriction map 
$S^3V \rightarrow H^0(M,\oo(3))$ has 
kernel of
dimension $2^g$.
Then there exists a unique invariant quartic $Q\in S^4_0 V$ with $M\subset \sing Q$.
\end{prop}

\pf
Since the restriction map is $J[2] = K \times \khat$-equivariant, its $K$-invariant part is 
$$
\ker\{
(S^3 V)^K \rightarrow H^0(M,\oo(3))^K\},
$$
which is therefore 1-dimensional, by \cite{M} proposition 3. By proposition \ref{partials}, therefore, we have a unique invariant quartic whose partial derivatives all vanish along $M$.
\qed

We shall apply this result to the images in $|2\Theta|$
of $\su(2)$ and the Kummer variety.

\begin{ex} $\bf g=2.$ \rm
Here $\phi$ is an isomorphism $\su(2) \idfy \bp^3$; see \cite{NR1}.
The Jacobian maps to the classical Kummer quartic surface, given by the 1-dimensional kernel of the surjective restriction map
$S^4_0 V\rightarrow H^0_0(J,8\Theta)$.

\end{ex}

\begin{ex} 
\label{exg=3}
$\bf g=3.$ \rm
Here $\dim S^3 V = 120$, $\dim H^0_+(J, 6\Theta) = 112$ and for nonhyperelliptic $C$---equivalently for ppav $J_C$ without vanishing theta-nulls---the restriction map $S^3 V \rightarrow H^0_+(J, 6\Theta)$ is surjective with 8-dimensional kernel.
So by proposition \ref{coble} there is a unique invariant quartic $Q_C\in S^4_0 V$ containing the Kummer variety in its singular locus---this is the original case of the proposition (see \cite{C} page 104).

On the other hand,  
Narasimhan and Ramanan \cite{NR2} showed that if $C$ is nonhyperelliptic then
$\phi:\su(2)\rightarrow |2\Theta|$ is an embedding whose image is an invariant  
quartic; in particular this quartic is singular along the Kummer and therefore coincides with $Q_C$.
It is called the {\it Coble quartic} of the curve.
\end{ex}

Alternatively, 
van Geemen and Previato \cite{vGP1} have shown that for curves without vanishing theta-nulls
$\phi(\su(2))$ is projectively normal in degree 4; in particular there is a surjective restriction map

\begin{equation}
\label{quarticrestriction}
S^4_0 V \rightarrow H^0_0(\su(2),\ll^4)
\end{equation}
where $H^0_0(\su(2),\ll^4)\subset H^0(\su(2),\ll^4)$ again denotes the subspace of $J_2$-invariants, with dimension (see \cite{OP})
$
\dim H^0_0(\su(2),\ll^4) = (3^g +1)/2
$. In the case $g=3$, therefore---where `no vanishing theta-nulls' means nonhyperelliptic---we see that $Q_C$ is the 1-dimensional kernel of~(\ref{quarticrestriction}).

\section{Cubic normality for genus 4}
\label{cubicnorm}

We shall now apply proposition \ref{coble} in the case of genus 4. Our main result is:

\begin{theo}
\label{cubicnormality}
For any curve $C$ of genus 4 without vanishing thetanulls the multiplication map
$
S^3 H^0(\su(2),\ll) \rightarrow H^0(\su(2), \ll^3)
$
is surjective.
\end{theo}

Using the Verlinde formula (see \cite{vGP1} or \cite{OW}) one observes that the respective dimensions of these spaces
are:
$$
\dim S^3 H^0(\su(2),\ll) = 816,
\qquad
\dim H^0(\su(2), \ll^3) = 800.
$$
(Note that the number $816 = 2^4 \times 51$ also comes from (\ref{dimS40}) and
proposition \ref{partials}.)
So from proposition \ref{coble} we deduce:

\begin{cor}
\label{specialq}
For any curve $C$ of genus 4 without vanishing theta-nulls
there exists a unique invariant quartic $Q_C \subset \bp V \dual \cong \bp^{15}$ with the property that $\phi(\su(2)) \subset \sing Q_C$.
\end{cor}

\begin{rems}\rm
\itemitem{(i)}
It is interesting to note what happens 
when the curve $C$ has a vanishing theta-null: assuming $C$ is nonhyperelliptic this vanishing theta-null is unique, i.e. $C$ has a unique semi-canonical pencil $g^1_3$, so that by \cite{B} proposition 2.6 the image $\phi(\su(2)) \subset |2\Theta| $ lies on a unique Heisenberg-invariant {\it quadric} $G$. In this case, therefore, the non-reduced quartic $Q_C = 2 G$ has the properties stated in the corollary; though we do not know that it is unique.

\itemitem{(ii)} When $C$ has no vanishing theta-nulls the same result \cite{B} proposition 2.6 tells us that $\phi(\su(2)) \subset |2\Theta| $ does not lie on any quadric, and it follows easily from this that the quartic in corollary \ref{specialq} is irreducible.
\end{rems}

Before proving theorem \ref{cubicnormality} we shall need some notation. 
Fix a theta structure and maximal level subgroups $K_g,\khat_g$ as in section \ref{invariantqs}.
A 2-torsion point $\eta \in K_g\subset J[2]$ acts linearly on $V$, and we shall denote by $V_{\eta}$ (resp. $V_{\eta}^-$) the $+1$- (resp. $-1$-)eigenspace. Since $K_g$ is an isotropic subgroup for the skew-symmetric Weil form on $J[2]$, and since the linear actions of two orthogonal 2-torsion points commute, the restriction map to the eigenspace $V_{\eta}$ maps $K_g$-invariant cubics to $K_g$-invariant cubics. Given any $\eta \in K_g$ we can choose an isomorphism
$$
K_g/\<\eta \> \idfy K_{g-1},
$$
and hence obtain a linear map
$$
\res_{\eta} : (S^3 V)^{K_g} \rightarrow (S^3 V_{\eta})^{K_{g-1}}.
$$
Furthermore, by the general theory of Prym varieties, once we choose a theta structure on the associated Prym $P_{\eta}$ the space $V_{\eta}$ becomes an irreducible $H_{g-1}$-module to which we may apply proposition \ref{partials}. 

\begin{lemm}
\label{injlem}
For $g\geq 3$ the map obtained by restriction to all the eigenspaces $V_{\eta}$ is injective:
$$
\sum \res_{\eta} : (S^3 V)^{K_g} \hookrightarrow \bigoplus_{\eta \in K_g\backslash\{0\}} (S^3 V_{\eta})^{K_{g-1}}.
$$
\end{lemm}

\begin{rem}\rm
Although we will not need the fact, one may note that the restriction map of $K_g$-invariant cubics to the $-1$-eigenspace $V_{\eta}^-$ is the zero map.
\end{rem}

\pf
Let $\eta = (1,0,\chi) \in 
K_g\backslash\{0\}$. By definition, the eigenspace $V_{\eta}$ is spanned by the vectors $X_{\sigma}$ such that $\chi(\sigma) =1$. These $\sigma$s form a subgroup of $\khat_g$ which is isomorphic to $\khat_{g-1}$; and $\{X_{\sigma}|\chi(\sigma) =1\}$ is a canonical basis for the action of $H_{g-1}$ on $V_{\eta}$. This means that the elements 
$$
\begin{array}{rclcl}
 X_{0}^3&&\\
 X_{0}^2 X_{0+\lambda}&&{\rm with}\ \chi(\lambda)=1,\\ 
 X_{\lambda} X_{\mu} X_{\lambda+\mu}&&{\rm with}\ \chi(\lambda)=\chi(\mu)=1,\\
\end{array}
$$
of $(S^3 V)^{K_{g}}$ map bijectively to a basis of $(S^3 V_{\eta})^{K_{g-1}}$. To prove the lemma it is therefore sufficient to observe that for any $\lambda,\mu \in \khat_g \backslash \{0\}$ there exists $\chi \in K_g \backslash \{0\}$ such that $\chi(\lambda) = \chi(\mu) = 1$. For $g\geq 3$ this is obvious.
\qed

{\it Proof of theorem \ref{cubicnormality}.} 
Choose a theta structure, so that the Heisenberg group $H_4$ acts on both spaces. Since the multiplication map is Heisenberg-equivariant it is enough to show surjectivity for $K$-invariant elements:
$$
m^K: (S^3 V)^K \rightarrow H^0(\su(2), \ll^3)^K.
$$
We have already observed that these spaces have dimensions 51 and 50 respectively, and so we have to show that the kernel is 1-dimensional. 

Consider a cubic $F\in \ker m^K$, and restrict it to the eigenspace $\bp V_{\eta}\dual$. Since the intersection of this eigenspace with $\phi(\su(2))$ is the Kummer image of the Prym variety $P_{\eta}$ (see \cite{NR3}), the restricted cubic $\res_{\eta}(F)$ is an element of
$$
U_{\eta} = \ker \{ (S^3 V_{\eta})^{K_3} \rightarrow H^0_+(P_{\eta}, 6\Xi)\}.
$$
By hypothesis $J_C$ has no vanishing theta-nulls, and it follows from the Schottky-Jung identities that $\peta$ has no vanishing theta-nulls.
By example \ref{exg=3}, therefore, $\dim U_{\eta} =1$.
It is therefore sufficient to show that for any $\eta$ the map 
$$
\res_{\eta}|_{\ker m^K } : \ker m^K \rightarrow U_{\eta}
$$
is injective. 

Suppose that $F\in \ker m^K \cap \ker \res_{\eta}$. Choose a nonzero $\zeta \in J[2]$ orthogonal to $\eta$ with respect to the Weil pairing; we shall show that $F\in \ker \res_{\zeta}$. Weil orthogonality 
implies that the intersection $V_{\eta,\zeta} =  V_{\eta}\cap V_{\zeta}$ is 4-dimensional (see, for example, \cite{vGP1}). 
By hypothesis $F$ vanishes in $V_{\eta}$ and hence in $V_{\eta,\zeta}$. On the other hand, $\res_{\zeta}(F) \in U_{\zeta}$; if this element is nonzero then it spans $U_{\zeta}$, and in particular its $\khat_3$-orbit spans all cubics in $\bp V_{\zeta}\dual$ vanishing on $\kum(P_{\zeta})$. We conclude that if  
$\res_{\zeta}(F) \not= 0$ then the singular locus of the Coble quartic of $P_{\zeta}$ contains $\bp V_{\eta,\zeta}\dual \cong \bp^3$, a contradiction.

By repeating this argument we deduce that if $F\in \ker m^K \cap \ker \res_{\eta}$ for {\it some} nonzero $\eta \in J[2]$ then the same is true for {\it any} $\eta \in J[2]$; by lemma \ref{injlem} this implies that $F=0$. 
\qed

{\it Proof of theorem \ref{theorem2}(1).}
Under the isomorphism of proposition \ref{partials}(1) the quartic $Q_C$
corresponds to a cubic $F$ whose restriction to each $\bp V_{\eta}\dual$ is nonzero, by the proof of theorem \ref{cubicnormality}. 
Therefore
since $K_g /\<\eta\> \cong K_{g-1}$ we see that $Q_C$ restricts to a nonzero invariant quartic in $\bp V_{\eta}\dual \cong \bp^7$ which is singular along $\kum (\peta) \subset \phi(\su(2))$. 

By example \ref{exg=3} (and since $\peta$ has no vanishing theta-nulls) this restriction is just the Coble quartic.
\qed

\section{Lines on the Coble quartic}
\label{cobleqs}

\begin{prop}
\label{seconcob}
Let $X$ be a nonhyperelliptic curve of genus 3. 
For $a,b \in J(X)$ the secant line $\overline{ab} \subset \bp^7$ of the Kummer variety lies on the Coble quartic $\sux(2)\subset \bp^7$
if and only if $a+b \in X-X$ or $a-b \in X-X$.
\end{prop}

\pf
We make use of the ruling of $\sux(2)\subset \bp^7$ by 3-planes (see \cite{OPP} or \cite{NR2}). First of all, if $a-b \in X-X$ then $a(p) = b(q) = x\in J^1(X)$ for some points $p,q,\in X$, and then the line $\overline{ab} \subset \bp^7$ is
precisely
the secant line $\overline{pq}$ of the image of the curve in the extension space $\bp(x) \subset \sux(2)$ (see \cite{OPP} for notation). If $a+b \in X-X$ the argument is similar.

So now suppose that $\overline{ab} \subset \sux(2)$. If $\overline{ab} \subset \bp(x)$ for some $x\in J^1(X)$, then (since $\bp(x)$ meets the Kummer in the image of $X\hookrightarrow \bp(x)$) $\overline{ab}$ is a secant line of the image curve, and one infers easily that $a+b \in X-X$ or $a-b \in X-X$.

We may therefore suppose that $\overline{ab}$ does not lie on any such 3-plane $\bp(x)$. Pick any stable bundle $E$ on the line $\overline{ab}$; then (since $g(X) = 3$) $E \in \bp(x)$ for some $x\in J^1(X)$. Under the embedding $\phi: \sux(2) \hookrightarrow |2\Theta|$ the linear subspace $\bp(x)$ is the linear system of divisors $D\in |2\Theta|$ containing the curve $x+ W_1 \cong X \subset J^2(X)$. We shall view the line $\overline{ab}$ as the pencil of $2\Theta$-divisors spanned by $\Theta_a + \Theta_{-a}$ and  $\Theta_b + \Theta_{-b}$,
and restrict it to $x+W_1$. 
Since there is a member of this pencil vanishing identically on the curve we see that there is an equality of effective divisors:

\begin{equation}
\label{divisor}
(\Theta_a + \Theta_{-a})|_{x+W_1} = (\Theta_b + \Theta_{-b})|_{x+W_1}.
\end{equation}
We note that $\Theta_a$ restricts to the line bundle $Kx^{-1}a$ on $x+W_1 \cong X$, {\it and moreover that $h^0(X,Kx^{-1}a) =1$}. To see this, observe that since $\deg Kx^{-1}a = 3$, $h^0(X,Kx^{-1}a) >1$ only if $Kx^{-1}a = K(-p) $ for some point $p\in X$. So $x = a(p)$ and we deduce that $x+W_1 \subset \Theta_a$. 
But then the Kummer image of $a$, i.e. the divisor $\Theta_a + \Theta_{-a}$, lies in the space $\bp(x)$ and therefore so does the line $\overline{ab}$, contrary to hypothesis.

We may thus write (\ref{divisor}) in the form:
$$
p_1+p_2+p_3 + q_1+q_2+q_3 = r_1+r_2+r_3 + s_1+s_2+s_3
$$
where $\{p_1+p_2+p_3\} = |Kx^{-1}a|$, $\{q_1+q_2+q_3\}= |Kx^{-1}a^{-1}|$ etc. If $p_1+p_2+p_3 = r_1+r_2+r_3$ or $s_1+s_2+s_3$ then we deduce that $a=\pm b$. 

Otherwise, we can find an equation of the form
$$
p_1+p_2+p_3 = r_1+r_2+ s_i
$$
for some $i=1,2,3$. But this means that 
$$
a-b = (p_1+p_2+p_3 ) - (r_1+r_2+ r_3) = s_i - r_3 \in X-X.
$$
Similarly, each such equation leads either to $a-b\in X-X$ or to $a+b\in X-X$.
\qed

\begin{cor} 
The secant lines of a 3-dimensional Kummer variety in $\bp^7$ cover $\bp^7$.
\end{cor}

\pf
The secant variety of the Kummer is irreducible and contains $\sux(2)$, since the latter is swept out by 3-planes meeting the Kummer in a curve. But the preceding proposition shows that the inclusion is proper; since $\sux(2) \subset \bp^7$ has codimension one, therefore, the result follows.
\qed

\section{Quadrisecant planes}
\label{quadrisecants}

In this section we shall recall the result of Beauville--Debarre \cite{BD}
which says that the Kummer variety of a Prym:
$$
\kum : P_{\eta} = {\rm Prym}(C,\eta) \rightarrow |2 \Xi | = \bp^{2^{g} -1}
$$
(where we shall take $C$ to have genus $g+1$)
possesses a 4-parameter family of quadrisecant $\bp^2$s, analogous to the trisecant lines of a Jacobian Kummer.

The base $\bb$ of the family is the fibre product:

\begin{equation}
\label{BDfamily}
\begin{array}{rcl}
\bb &\rightarrow & S^4 \ctil\\
&&\\
\downarrow &&\downarrow \hbox{Abel-Prym}\\
&&\\
P_{\eta}^- & \map{sq} & P_{\eta}\\
\end{array}
\end{equation}
where $sq$ is the squaring map $a\mapsto a^{2}$.
An element of $\bb$, in other words, is a pair $(a,\Gamma) \in P_{\eta}^-\times 
S^4 \ctil$ such that $a^{2} = [\Gamma]$ (see section \ref{pryms}).
For any points $p,\ldots,q\in \ctil$ and $a\in P_{\eta}^-$, let us write:
$$
\begin{array}{rcl}
\<p\>_a &=& \kum([p] -a)\in |2 \Xi | = \bp^{2^{g} -1},\\
\<p,\ldots ,q\>_a &=& \hbox{linear span of $\<p\>_a,\ldots ,\<q\>_a$}.\\
\end{array}
$$
The following fact then results from \cite{BD} propositions 1 and 2.

\begin{prop}
\label{linepairs}
For each $(a,p+q+r+s) \in \bb$, the four points 
$
\<p\>_a,\<q\>_a,\<r\>_a,\<s\>_a  \in \kum\ P_{\eta} \subset \bp^{2^{g} -1}
$
are coplanar.
\end{prop}

We now make the following observation which relates this family to the parameter space $\ff$ of Fay trisecants (see \cite{OPP}, section 2). First, $\ff \subset J^2(C)\times S^4 C$ consists of pairs $(\lambda, D)$ such that $\lambda ^2 = \oo(D) \in J^4(C)$; and we consider the 16:1 cover 
$$
\begin{array}{rcc}
\gg \cup \gg^- = \fftil &\subset & J^2(C) \times S^4 \ctil\\
&&\\
16:1 \downarrow &&\downarrow {\rm id}\times \nm_{\pi}\\
&&\\
\ff & \subset& J^2(C) \times S^4 C.\\
\end{array}
$$
There is a map $\fftil \rightarrow \peta \cup \peta^-$ defined by $(\lambda,\Gamma) \mapsto \pi^*\lambda^{-1} \otimes \oo(\Gamma)$, and replacing $\Gamma = p+q+r+s$ by $\sigma(p) +q+r+s$, i.e. transposing one point of $\Gamma$ by sheet-interchange over $C$, switches the image from $\peta$ to $\peta^-$ and vice versa. This shows that $\fftil$ has two connected components, and accordingly we have written $\fftil = \gg \cup \gg^-$. Thus, by definition
$$
\gg^- = \{ (\lambda ,\Gamma) \in J^2(C) \times S^4 \ctil \,|\, \lambda^2 = \nm_{\pi} \Gamma;\ \pi^*\lambda^{-1} \otimes \oo(\Gamma)\in \peta^-\}.
$$
If we write $\Gamma = p+q+r+s$ again then the sixteen points of $\peta \cup \peta^-$ coming from the fibre of $\gg \cup \gg^-$ over $(\lambda , \nm_{\pi} \Gamma) \in \ff$ are (modulo $\sigma$):

\begin{equation}
\label{FBDtable}
\begin{array}{|l|l|}
\hline
&\\
\qquad \peta^-&\qquad \peta\\
\hline
&\\
\pi^* \lambda^{-1} (p+q+r+s)=a  & \pi^* \lambda^{-1}(\sigma p+q+r+s) = [p] -a\\
\pi^* \lambda^{-1} (\sigma p+\sigma q+r+s) & \pi^* \lambda^{-1}( p+\sigma q+r+s) = [q] -a\\
\pi^* \lambda^{-1} (\sigma p+q+\sigma r+s) & \pi^* \lambda^{-1}( p+q+\sigma r+s) = [r] -a\\
\pi^* \lambda^{-1} (\sigma p+q+r+\sigma s)& \pi^* \lambda^{-1}( p+q+r+\sigma s) = [s] -a\\
&\\
\hline
\end{array}
\end{equation}
Note that we define $a = \pi^* \lambda^{-1} (\Gamma)$ here, and that $-a = \sigma a = \pi^* \lambda^{-1} (\sigma \Gamma)$. Observe also that $a^2 = [\Gamma]$ so that we have a map $g: \gg^- \rightarrow \bb$ sending $(\lambda, \Gamma) \mapsto (a ,\Gamma)$. In addition, we have a map 
$$
\begin{array}{rcl}
f: \ff &\rightarrow & S^4 \kum(\peta) \\
x & \mapsto & \{\hbox{fibre of $\gg$ over $x$}\}/\sigma.\\
\end{array}
$$
In other words, $f$ maps $x = (\lambda , \nm_{\pi} \Gamma)$ to the four points in the right-hand column of (\ref{FBDtable}).

In conclusion, the identifications in the right-hand column of the table show that the following diagram commutes:

\begin{equation}
\label{FBDdiagram}
\begin{array}{rcl}
\gg^- & \map{g} & \bb\\
&&\\
8:1 \downarrow &&\downarrow {\hbox{Beauville-Debarre}\atop\hbox{quadrisecant planes}}\\
&&\\
\ff & \map{f}& S^4 \kum(\peta) \hookrightarrow \gr_3 H^0(2\Xi).
\end{array}
\end{equation}

\begin{rems}\rm
\itemitem{(i)}
It is easy to see that $\gg^-$ is \'etale on both $\ff$ and $\bb$ (with degree 2 in the latter case); in particular both spaces have the same image in $S^4 \kum(\peta)$.
\itemitem{(ii)}
We shall make use of diagram (\ref{FBDdiagram})
in section \ref{prymG3}. The fact that $\ff$ parametrises both the trisecants of the Jacobian Kummer and the quadrisecants of the Prym Kummers
is rather striking, but will not play any role in this paper.
\end{rems}

We now
choose a $g^1_3$ on $C$. The Recillas correspondence then determines a tetragonal curve $(X,g^1_4)$ and identifies $J(X) \cong P_{\eta}(C)$. In particular, we can now view $\ctil \subset S^2 X$ as in diagram (\ref{recdiagram}); we shall denote this inclusion map by $p\mapsto p_1 + p_2$ for $p\in \ctil$.

Restricting to the case $g=3$, we shall view $X\subset \bp^2$ in its canonical embedding, and consider lines $\overline{p} = \overline{p_1 p_2} \subset \bp^2$, for each $p\in \ctil$.

\begin{prop}
Let $g=3$ and $\sux(2) \subset \bp^7$ be the Coble quartic with singular locus $\kum\ P_{\eta}(C)$. Then for $(a, p+q+r+s) \in \bb$ the following are equivalent:
\begin{enumerate}
\item
$\<p,q\>_a \subset \sux(2)$;
\item
$\<r,s\>_a \subset \sux(2)$;
\item
$\overline{p}\cap \overline{q} \in X\subset \bp^2$ or $\overline{r}\cap \overline{s} \in X\subset \bp^2$.
\end{enumerate}
\end{prop}

\pf
We consider the line $\<p,q\>_a$. By proposition \ref{seconcob} this lies on $\sux(2)$ if and only if the sum $[p]+[q]-2a$ or the difference $[p]-[q]$ of the points $[p]-a, [q]-a \in P_{\eta}(C) = J(X)$ lies on the surface $X-X$. In view of the diagram (\ref{recdiagram}) we have

\begin{equation}
\label{one}
[p]-[q] = p_1 +p_2 - q_1-q_2 \in J(X);
\end{equation}
whilst by definition $2a = [p+q+r+s] \in P_{\eta}$ is identified with the point $p_1+p_2 + \cdots + s_1+s_2 - 2K_X \in J(X)$, so that 

\begin{equation}
\label{two}
[p]+[q] - 2a = K_X -r_1-r_2 -s_1-s_2 \in J(X).
\end{equation}
By (\ref{one}), $[p]-[q]\in X-X$ if and only if 
$$
p_1+p_2 +v \sim q_1+q_2+ u
$$
for some $u,v \in X$. If these divisors are equal then $p_i = q_j$ for some $i,j\in \{1,2\}$; in which case $\overline{p} \cap \overline{q}$ lies on $X\subset \bp^2$ trivially. Otherwise the two divisors span a $g^1_3 = |K_X(-s)|$ for some $s\in X$; and so in this case $\overline{p} \cap \overline{q}$ is the point $s\in X$. Conversely, if $\overline{p} \cap \overline{q}\in X$ then the same argument shows that $[p]-[q] \in X-X$.

By (\ref{two}), on the other hand, $[p]+[q] -2a \in X-X$ if and only if 
$$
r'_1+r'_2 +v \sim s_1+s_2+ u
$$
for $u,v \in X$, where $r_1+r_2+r'_1+r'_2 \in |K_X|$. And now, since $\overline{r'_1 r'_2} = \overline{r_1 r_2} = \overline{r}$, the same argument as before shows that $[p]+[q] -2a \in X-X$ if and only if $\overline{r} \cap \overline{s}\in X$.

We conclude that parts 1 and 3 are equivalent to each other and hence also to part 2.
\qed

From the preceding proposition we can deduce necessary and sufficient conditions for the quadrisecant plane $\<p,q,r,s\>_a \subset \bp^7$ to lie on $\sux(2)$: since $\sux(2)\subset \bp^7$ has degree 4 this is equivalent to $\sux(2)$ containing all of the six lines $\<p,q\>_a, \ldots, \<r,s\>_a$. Let 
$$
\Delta(p) = \overline{q} \cup \overline{r} \cup\overline{s}\subset \bp^2,
$$
with $\Delta(q),\Delta(r),\Delta(s)$ defined similarly. Then from proposition \ref{linepairs} we see that $\sux(2)$ contains the six lines if and only if {\it either} one of the triangles $\Delta(p),\ldots , \Delta(s)$ has its vertices on $X\subset \bp^2$ {\it or} one of the lines, $\overline{p}$ say, meets the remaining three in points of $X\subset \bp^2$. 

{\it To summarise:} $\<p,q,r,s\>_a \subset \sux(2)$ if and only if the four lines $\overline{p},\overline{q},\overline{r},\overline{s} \subset \bp^2$ form a configuration with respect to the canonical curve of one of the two forms below:

$$
\epsfxsize=1.5in\epsfbox{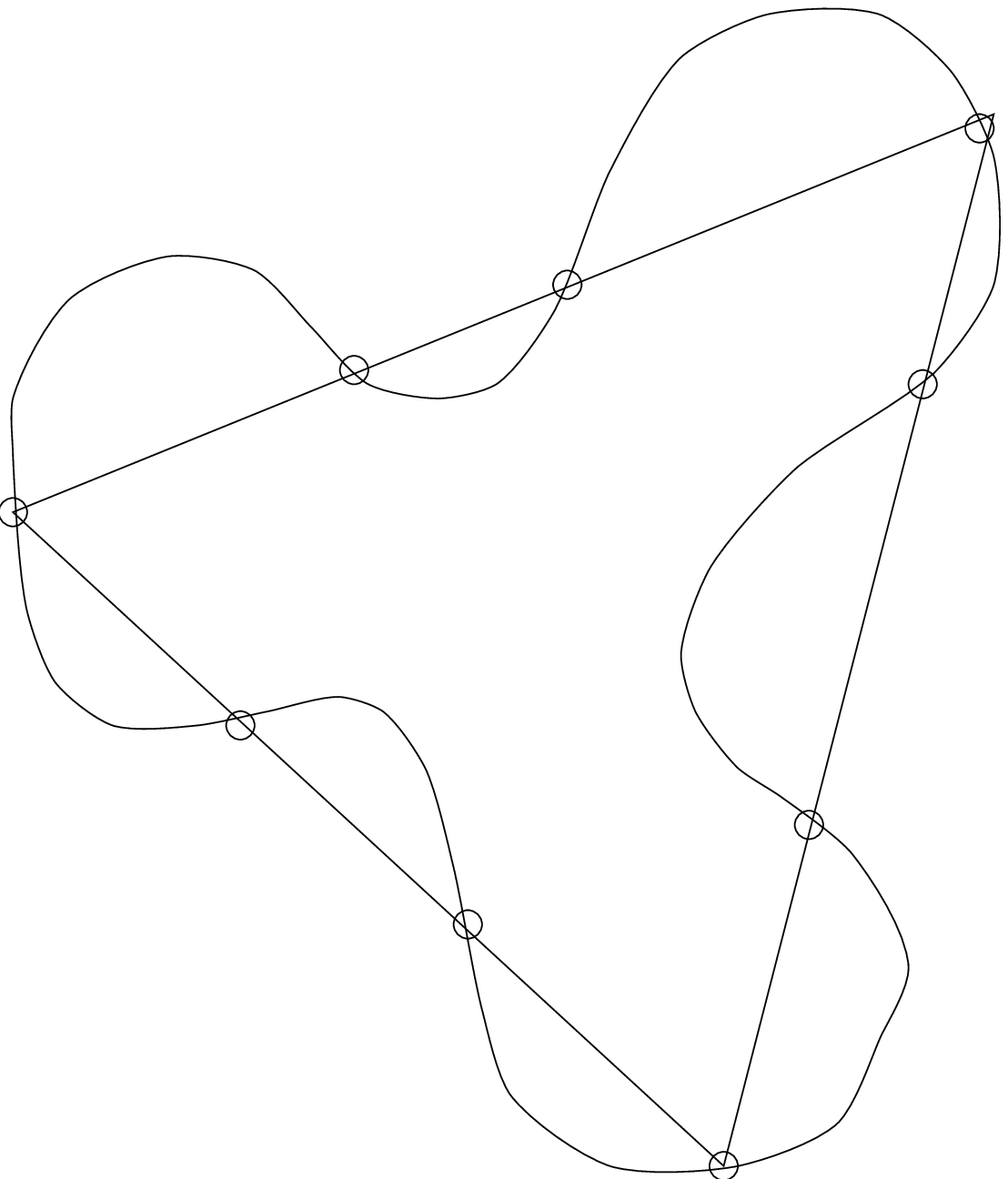}
\qquad\qquad\qquad
\epsfxsize=2in\epsfbox{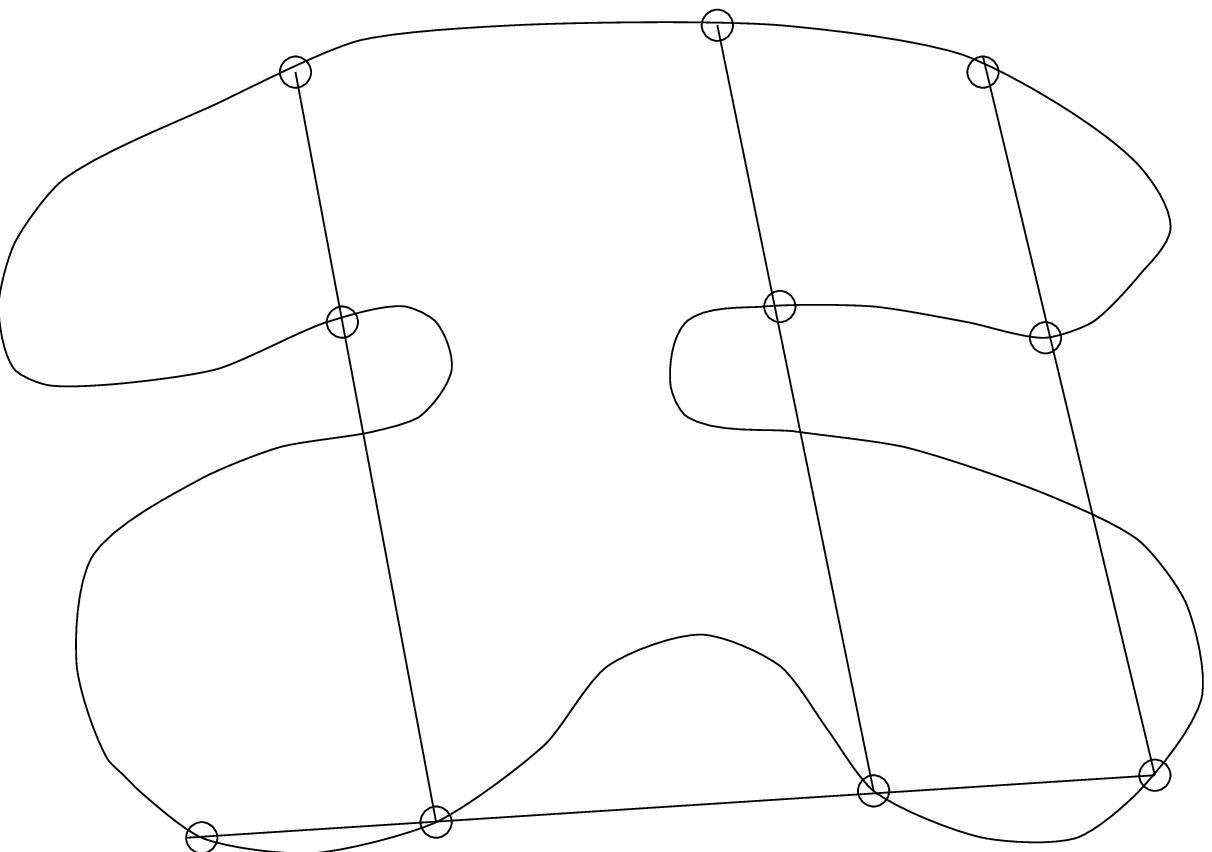}
$$

\begin{cor}
The generic quadrisecant plane $\<p,q,r,s\>_a$ is not contained in $\sux(2)$. 
\end{cor}

\pf
By the Recillas construction $X$ comes equipped with a $g^1_4$, and each of the divisors $p_1 +p_2, \ldots , s_1 + s_2$ is contained in a divisor of this pencil. If the $g^1_4$ is canonical then it is obtained by projection of $X$ away from a point of the plane off the curve; in this case we see that the configurations above can never occur, i.e. that {\it no} quadrisecant plane $\<p,q,r,s\>_a$ is contained in $\sux(2)$.

More generally, we can count parameters. There can be at most a 2-parameter family of triangles on the left: one parameter for any side, and one for the opposite vertex. There is a 1-parameter choice of the remaining line; and hence at most a 3-dimensional family of quadrisecants can give rise to the configuration on the left.
There is also at most a 3-dimensional family of quadrisecants giving rise to the right-hand configuration. 
\qed

\section{Abel-Jacobi stratification of $|2\Theta|$}
\label{abeljacobi}

Let $W_d \subset J^d$ be the subscheme of special line bundles of degree $d$, $0\leq d \leq g-1$, i.e. the Abel-Jacobi image of $S^d C$. 
We introduce the following sets of divisors on the Jacobian:

\begin{equation}
\label{ajsets}
\begin{array}{rcl}
G_{d+1} &=& \{\ D\ | \ x+W_{g-1-d} \subset {\rm supp}\ D\ \hbox{for some $x\in J^d$}\ \} \subset |2\Theta|,\\
&&\\
G_{g-d}\dual &=& \{\ D\ | \ -x+W_d \subset {\rm supp}\ D\ \hbox{for some $x\in J^d$}\ \} \subset |\ll|.\\
\end{array}
\end{equation}
It will often be convenient to identify $|\ll| \idfy |2\Theta|$ via a choice of theta characteristic as in diagram (\ref{thetachoice}) of section \ref{prelims}. It is easy to see that this induces identifications $G_i\dual \idfy G_i$ for each $i=1,\ldots, g$; note also that this is unaffected by the choice of theta characteristic as all the sets (\ref{ajsets}) are $J[2]$-invariant.

Obviously $G_g = |2\Theta|$ and $G_0$ is empty, whilst by the theorem of the square $G_1 = \kum(J)$ where $\kum$ is the Kummer map
$$
\begin{array}{rcl}
\kum:J &\rightarrow & |2\Theta|\\
  x &\mapsto  & \Theta_x + \Theta_{-x}.\\
\end{array}
$$
We thus have (isomorphic) filtrations:

\begin{equation}
\label{AJstrat}
\begin{array}{rcl}
\kum(J) &= G_1 \subset G_2 \subset \cdots \subset G_{g-1}& \subset  |2\Theta|,\\
&&\\
\kum(J^{g-1}) &= G_1\dual \subset G_2\dual \subset \cdots \subset G_{g-1}\dual& \subset |\ll| = |2\Theta|\dual.\\
\end{array}
\end{equation}
To interpret these sets scheme-theoretically we consider maps
$$
\begin{array}{ccc}
J^d \times S^d C&\map{\alpha_d}&J\\
&&\\
\downarrow \pi&&\\
&&\\
J^d&&\\
\end{array}
$$
where $\alpha_d: (x,D) \mapsto x^{-1} \otimes \oo (D)$ and $\pi$ is  projection to the first factor; and set
$Q_d = \pi_* \alpha_d^* \ll.$
In the notation of the appendix and proposition \ref{alphabeta}, 
$$
Q_d(x) = H^0(S^d C, \ll_x) = H^0(\bp (x), \iic^{d-1}(d))
$$ 
at each $x\in J^d$. 

{\it In the remainder of this section we shall assume that $g\leq 4$.} More generally (see proposition \ref{conjOK}) we shall assume the validity of conjecture \ref{400}.)

It then follows that $Q_d$ is locally free and the restriction map $\alpha_d^*:\oo_{J^d} \otimes H^0(J,\ll)\rightarrow Q_d$ is surjective. 
We therefore obtain dual short exact sequences of vector bundles on $J^d$ (in which $N_d = \ker \alpha_d^*$ by definition):

\begin{equation}
\label{NQsequence}
\begin{array}{rcl}
0\rightarrow N_d \rightarrow &\oo_{J^d} \otimes H^0(J,\ll)& \map{\displaystyle
\alpha_d^*} Q_d
 \rightarrow 0\\
&&\\
0\rightarrow Q_d\dual \rightarrow &\oo_{J^d} \otimes H^0(J^{g-1},2\Theta)& \rightarrow N_d\dual
 \rightarrow 0.\\
\end{array}
\end{equation}
{\it We define $G_{g-d}\dual$ to be the ruled variety $\im \{\, \bp N_d \subset J^d \times |\ll| \rightarrow |\ll|\, \}$, and
$G_{d+1}\subset |2\Theta|$ to be the ruled variety $\im \{\, \bp Q_d\dual\subset J^d \times |2\Theta|\rightarrow |2\Theta|\, \}$.}

Thus, for $x\in J^d$ the subset $\bp N_d(x)\subset |\ll|$ consists of divisors containing $-x+W_d$. In a moment we shall verify the less obvious fact that 

\begin{equation}
\label{Qdxdual}
\bp Q_d (x)\dual = \{\ D\in |2\Theta|\ | \ x+W_{g-1-d} \subset {\rm supp}\ D\ \ \}
\end{equation}
so that this definition is compatible with (\ref{ajsets}).

\begin{rem}\rm
\label{dualdims}
By \ref{400}(1) we have
$\rk Q_d  = \sum_{i=0}^d {g\choose i}$.
It follows from this and the fact that $h^0(J,\ll) = 2^g$
that 
$$
\rk N_d = \sum_{i=d+1}^g {g\choose i} = \rk Q_{g-1-d}.
$$
Note in particular that $\dim \bp N_1 = \dim \bp Q_{g-2}\dual = 2^g -2$, and we can therefore expect that $G_{g-1}\subset |2\Theta|$ is a hypersurface. We shall be particularly interested in this case in what follows.
\end{rem}

Recall also that a choice of $\kappa \in \vartheta (C)$ gives rise to 
an isomorphism
$|2\Theta| \cong |2\Theta|\dual = |\ll|$ (see section \ref{prelims}).
Given this choice we can prove the following 
polarity relation:

\begin{prop}
\label{perp}
Fix $x\in J^d$ and 
theta characteristic $\kappa \in \vartheta (C)$. Then
the induced isomorphism $H^0(J^{g-1}, 2\Theta) \idfy H^0(J,\ll) =H^0(J^{g-1}, 2\Theta) \dual $ (up to scalar)
restricts to
$$
Q_d(x) \dual\idfy N_{g-1-d}(\kappa x^{-1}) = N_d(x)^{\perp}. 
$$
\end{prop}

\pf
By definition, the annihilator of $\bp N_d(x)\subset |\ll|$ in $|\ll|\dual = |2\Theta|$ is the linear span of 
$
\kum(W_d - x)\subset |2\Theta|.
$
Likewise $\bp N_{g-1-d}(\kappa x^{-1})$ has annihilator spanned by 
$
\kum(W_{g-1-d} -\kappa +x) \subset |2\Theta|.
$
Now, under the map $m:J\times J \rightarrow J\times J$, $(u,v)\mapsto (u+v,u-v)$ of section \ref{prelims} we note that 
$$
(W_d - x)\times (W_{g-1-d} -\kappa +x) \map{m} \Theta_{\kappa} \times J
$$
and hence by (\ref{riemannxi}) the image is killed by the form $\xi_{\kappa}$ defining the isomorphism $|\ll| \cong |\ll|\dual$. It follows that 
$$
N_d(x)\perp N_{g-1-d}(\kappa x^{-1})
$$
under this pairing. But from remark \ref{dualdims} we see that 
$$
\dim N_d(x) + \dim N_{g-1-d}(\kappa x^{-1}) = 2^g,
$$
and so $N_{g-1-d}(\kappa x^{-1}) = N_d(x)^{\perp}$.

On the other hand,
let $V=H^0(J,\ll)$ as usual. Since $N_d = N_{g-1-d}^{\perp}$ there is an isomorphism
$Q_d = (\oo_J\otimes V) /N_d \idfy (\oo_J\otimes V) /N_{g-1-d}^{\perp} = N_{g-1-d}\dual$.
\qed

In the following corollary, let $V_{\eta}^{\pm}\subset V$ be the eigenspaces of the action of $\eta \in J[2]$ on $V$ (or more precisely a lift of $\eta$ to the Heisenberg group); and note that under the pairing induced by any theta characteristic as above, $V_{\eta}^+ \perp V_{\eta}^-$. It is then easy to check the following fact using remark \ref{dualdims} and proposition \ref{perp}.

\begin{cor}
\label{intdims}
For each $x\in J^d$ and $\kappa \in \vartheta(C)$, $\eta\in J[2]\backslash \{0\}$, we have
$$
\dim(N_d(x) \cap V_{\eta}^+) - \dim(N_{g-1-d}(\kappa x^{-1}) \cap V_{\eta}^-)
= 2^{g-1} - \sum_{i=0}^d {g\choose i}.
$$
\end{cor}


The first step in the direction of calculating the degrees of the varieties $G_{d+1} \subset |2\Theta|$ is the observation that the self-intersection of the hyperplane class in $\bp N_d$ is the top Segre class of the bundle, and that in view of the exact sequence (\ref{NQsequence}) (since the central term is trivial) this in turn is the top Chern class of $Q_d$:
$$
\oo_{\bp N_d}(1)^{\dim \bp N_d} = s_g(N_d) = c_g(Q_d) \in H^{2g}(J^d,\z) \cong \z.
$$
In other words, we have

\begin{equation}
\deg G_{g-d}\  \times \ \Bigl( 
{\hbox{no. of translates of $W_d$} \atop \hbox{in a general $D\in G_{g-d}$}}
\Bigr)
=
c_g(Q_d) \in \z.
\end{equation}
We shall restrict attention to the case $d=1$; the Chern classes of $Q_1$ are easily computed by comparison with the direct image sheaf
$$
Q' = \pi_* (\pp^2 \otimes pr_C^*K)
$$
where $\pp \rightarrow J^1 \times C$ is a Poincar\'e bundle. By proposition \ref{100}(1) we have 

\begin{equation}
\label{QNQ'}
Q_1 = \nn \otimes Q'
\end{equation}
for some line bundle $\nn \in \pic(J^1)$. To determine $\nn$, note that if we fix a base-point $p_0 \in C$ the Poincar\'e bundle $\pp$ is uniquely determined by requiring that $\pp|_{J^1 \times \{p_0\}}$ be trivial. This will be satisfied, with

\begin{equation}
\label{computingN1}
\pp^2 \otimes pr_C^*K = \pi^*\nn^{-1} \otimes \alpha_1 ^* \ll,
\end{equation}
if $\nn$ is chosen so that $\pi^*\nn^{-1} \otimes \alpha_1 ^* \ll$ is trivial on $J^1 \times \{p_0\}$, i.e. if we let

\begin{equation}
\label{computingN2}
\nn = \alpha_{p_0}^* \ll = t_{-p_0}^* \ll
\end{equation}
where $t_{-p_0} : J^1 \rightarrow J$ denotes translation by $p_0$. Applying $\pi_*$ to both sides of (\ref{computingN1}) we obtain (\ref{QNQ'}).

The Chern character of $Q'$ can be computed by Grothendieck-Riemann-Roch (see \cite{ACGH} chapter VIII) and is $g+1 - 4\theta$, where $\theta$ denotes the fundamental class of $\Theta$ in the Jacobian. We therefore conclude from (\ref{QNQ'}) and (\ref{computingN2}) that 

\begin{equation}
\label{chernQ}
ch(Q_1) = (g+1-4\theta) e^{2\theta}.
\end{equation}
Expanding the exponential and applying Newton's formula (\cite{F} page 56) we deduce:

\begin{prop}
The Chern classes $c_n = c_n(Q_1)$ are given recursively by 
$$
c_n ={1\over n} (c_{n-1}p_1 - c_{n-2} p_2 + \cdots + (-1)^{n-2} c_1 p_{n-1} + (-1)^{n-1}p_n)
$$
where $p_n = (g+1 - 2n) (2\theta)^n$.
\end{prop}

Although it is not easy to obtain a closed formula for $c_g(Q_1)$ it is readily calculated on a computer, using the above proposition, and the first few values are:

$$
\begin{array}{r|l}
g & c_g(Q_1) \\
&\\
\hline
&\\
3&32\\
4&384\\
5&4096\\
6&56320\\
7&872448\\
8&15368192\\
\end{array}
$$
We conjecture that at least for $g=4$, $\deg G_{g-1} = 4$. One might hope to use theorem \ref{theorem2}(2) to verify this.

\section{Segre stratification of $\su(2)$}
\label{segrestrat}

We shall show in this section that the Abel-Jacobi stratification (\ref{AJstrat}) induces, via $\phi$, the {\it Segre stratification} of $\su(2)$, i.e. the stratification by maximal degree of a line subbundle. (See \cite{LN}.) For any vector bundle $E$ on $C$ of rank 2 and degree 0
let
$$
n(E) = {\rm min}\{\, n\in \z\, | \,\hbox{$\exists$ line subbundle $L\subset E$ with $\deg L = -n$}\,\}.
$$
This function is nonnegative on semistable bundles; it
is also lower semicontinuous on families and determines a filtration of moduli space by closed subvarieties:

\begin{equation}
\label{Segrestrat}
J/{\pm} \cong \su(2)_0 \subset \su(2)_1 \subset \cdots \subset \su(2)_{[g/2]} = \su(2)
\end{equation}
where $\su(2)_d = \{E\in \su(2)|n(E) \leq d\}$. The right-hand equality in (\ref{Segrestrat}) follows from Nagata's theorem \cite{N}, \cite{L}: every ruled surface of genus $g$ admits a section with self-intersection at most $g$.

\begin{prop}
\label{strattheorem}
For $E\in \su(2)$, 
$n(E)\leq d$
if and only if
$\phi(E) \in G_{d+1}$.
In other words, $\su(2)_d = \phi^{-1} G_{d+1}$ for $0\leq d \leq [g/2]$.
\end{prop}

By definition $n(E)\leq d$
if and only if there exists $\xi \in J^d$ such that $h^0(C,\xi \otimes E)>0$; and that this is equivalent to $\phi(E) \in G_{d+1}$
can be restated in the following form, proved in \cite{OPP}, section 4.
For any $E\in \su(2)$ and $\xi \in J^d$:

\begin{equation}
\label{killingD}
h^0(C,\xi \otimes E)>0
\quad
\Longleftrightarrow
\quad
{h^0(C,\xi(D) \otimes E)>0\atop
\forall\,\, D\in S^{g-1-d}C.}
\end{equation}

Note that the case $d=0$ of proposition \ref{strattheorem}
just says that $\kum(J)\subset |2\Theta|$ comes from the semistable boundary of $\su(2)$; at the other end we see that
$$
\phi(\su(2)) \subset G_{[g/2]+1}.
$$

\begin{rems}\rm
\itemitem{(i)}
In a moment (see corollary \ref{G2}) we shall show that $G_2 \subset \phi(\su(2))$. By remark \ref{dualdims} $G_2$ is ruled by $\bp^g$s; this is the `$g$-plane ruling' of \cite{OPP} \S1.
\itemitem{(ii)}
For $g=2$ we know $\phi(\su(2)) = G_2 = |2\Theta|$; while for each $g\geq 3$ we have $\phi(\su(2)) \subset G_{g-1}$. In the case $g=3$ equality holds and $G_2\subset \bp^7$ is the Coble quartic of the curve.
\end{rems}

By Serre duality the projective space $\bp (x) := |Kx^2|\dual$, for $x\in J^d$, can be identified with the space $\bp H^1(C, x^{-2})$ of nontrivial extensions
$$
\ses{x^{-1}}{E}{x}
$$
up to isomorphism of $E$. It therefore has a rational map to $\su(2)$ (we shall see in a moment that the generic extension is semistable) which has been described in detail by Bertram, Lange--Narasimhan and others \cite{Bert}, \cite{LN}. 

\medskip\noindent
{\it Notation:}
For each $x\in J^d$ we shall denote by $\eee : \bp (x) \rightarrow \su(2)$ the rational moduli map, and write $E=\eee(e)$ for the bundle corresponding to 
a point $e\in \bp (x)$.
We shall denote by $\secant^n C \subset \bp(x)$ the variety of $n$-secants $\dbar \cong \bp ^{n-1}$, for $D\in S^n C$, of $\lambda_{|K x^2|}:C\hookrightarrow \bp(x)$. 
Finally, in the lemma below we shall consider the blow-up of $\bp(x)$ along these secant varieties: we shall denote by $\sectil^n C$ the proper transform of any $\secant^n C$ in a blow-up with lower dimensional centre, and by $\ss_n$ the exceptional divisor of the blow-up along $\secant^n C$ or $\sectil^n C$.
 
\medskip

It is easy to show (see for example \cite{LN} proposition 1.1) that 

\begin{equation}
\label{dbar}
e\in \dbar
\quad
\Longleftrightarrow
\quad
x(-D) \subset E;
\end{equation}
and in particular that for all $|n| \leq d-1$

\begin{equation}
\label{seccondition}
n(E)\leq n
\quad
\Longleftrightarrow
\quad
e\in \secant^{d+n}C 
\end{equation}
Note that the requirement $n\leq d-1$ arises here because $n(E) \leq d$ for all extension classes in $\bp (x)$ since $x^{-1}$ is always a line subbundle.

Applying (\ref{dbar}) when $\deg D =d$ says that $\secant^d C$ maps---away from $\secant^{d-1} C$---to the semistable boundary $J/\pm$, each $d$-secant plane $\dbar$ contracting down to the S-equivalence class of $E= x(-D) \oplus x^{-1}(D)$.

Applying (\ref{dbar}) when $\deg D <d$---and note that by (\ref{seccondition}) $\secant^{d-1}C$ is the exceptional locus of the rational map $\eee$---we see that in this case $\dbar$ transforms to the image $\eee \bp (x(-D))$.

The picture is clarified by the following results of Bertram \cite{Bert}.

\begin{lemm}
\label{bertram}
For $x\in J^d$ the rational map $\eee : \bp (x) \rightarrow \su(2)$ has the following properties. 
\begin{enumerate}
\item
$\eee^*: H^0(\su(2),\ll) \rightarrow H^0(\bp(x),\iic^{d-1}(d))$;
\item
$\eee$ resolves to a morphism $\eeetil$ of the $(d-1)$-st blow-up:
$$
\begin{array}{r}
\bp (x) \leftarrow \bl_C \leftarrow \bl_{\sectil^2 C} \leftarrow \cdots \leftarrow \bl_{\sectil^{d-1}C} = \widetilde\bp (x)\\
\\
\downarrow \eeetil\\
\\
\su(2)\\
\end{array}
$$
\item 
For $1\leq n\leq d-1$ we have
$\displaystyle \eeetil(\ss_n) = \bigcup_{D\in S^n C} \eeetil\, \widetilde\bp (x(-D))$.
\end{enumerate}
\end{lemm}

Combining these observations with proposition \ref{strattheorem}, the situation can be further summarised in the following diagram:

$$
\begin{array}{rcccccl}
\sectil^d C&\subset \cdots\subset &\sectil^{2d-1} C&\subset &\widetilde \bp (x)&&\\
&&&&&&\\
\downarrow&&\downarrow&&\downarrow \eeetil&&\\
&&&&&&\\
J/\pm  \cong \su(2)_0&\subset \cdots  \subset &\su(2)_{d-1}& \subset& \su(2)_d& \subset \cdots \subset &\su(2)\\
&&&&&& \\
 \downarrow&&\downarrow && \downarrow &&\downarrow \phi\\
&&&&&&\\
\kum(J) = G_1 &\subset  \cdots \subset& G_d & \subset& G_{d+1}&\subset \cdots \subset& G_{[g/2]+1}\subset |2\Theta|\\
\end{array}
$$

\medskip
\begin{rem} \rm
Since $\dim \bp (x) = g + 2d -2$ we see that $\secant^{d+n} C = \bp (x)$ as soon as $n\geq (g-1)/2$, i.e. $n\geq [g/2]$. This means that for such $n$ the top row of the diagram terminates, which is consistent with the fact that $\su(2) = \su(2)_{[g/2]}$.
\end{rem}

Now consider the linear pull-back map $\eee^* : |\ll| \rightarrow |\iic^{d-1}(d)|$ and its dual $(\eee^*)\dual: |\iic^{d-1}(d)|\dual\rightarrow
|2\Theta|$. We will check that it fits into the following commutative diagram, where $\alpha_x$ and $\beta_x$ are as defined in the appendix ((\ref{alphax}) and (\ref{betax}) respectively):

\begin{equation}
\label{abdiagram}
\begin{array}{ccc}
S^d C& \map{\displaystyle\beta_x} & |\iic^{d-1}(d)|\dual \\
&&\\
\alpha_x\downarrow& &\downarrow(\eee^*)\dual\\
&&\\
J& \map{\displaystyle\kum}& |2\Theta|\\
\end{array}
\end{equation}

To see this, first note that $D\in S^d C$ maps under $\kum \circ \alpha_x$ to the split divisor $\Theta_{x(-D)}+\Theta_{x^{-1}(D)}$, i.e. to the image of the semistable equivalence class of $x(-D)\oplus x^{-1}(D)$ under the morphism $\phi:\su(2) \rightarrow |2\Theta|$. 

On the other hand, $\beta_x(D)$ can be viewed as the set of divisors in $|\iic^{d-1}(d)|$ containing the linear span $\dbar \subset \bp(x)$. But we have seen above that this linear span contracts under $\eee$ to the point $x(-D)\oplus x^{-1}(D)$; thus $(\eee^*)\dual\circ \beta_x(D)$ coincides with $\phi(x(-D)\oplus x^{-1}(D))$, and the diagram commutes.

\begin{cor}
\label{completeness}
For $x\in J^d$, suppose that $\alpha_x^*$ is surjective (see \ref{400}(2)).
Then
the pull-back $\eee^* :H^0(\su(2),\ll) \rightarrow H^0(\bp(x), \iic^{d-1}(d))$ is surjective; in other words, the rational map $\eee: \bp(x) \rightarrow \su(2) \hookrightarrow |2\Theta|$ cuts out the complete linear system $|\iic^{d-1}(d)|$.
\end{cor}

\pf
In view of the identification $H^0(\su(2),\ll) = H^0(J,\ll)$, the commutativity of (\ref{abdiagram}) implies commutativity of
$$
\begin{array}{rcl}
H^0(J,\ll) & \map{\displaystyle
\eee^*} & H^0(\bp(x), \iic^{d-1}(d))\\
&&\\
&\alpha_x^* \searrow& \downarrow \beta_x^*\\
&&\\
&&H^0(S^d C, \ll_x).\\
\end{array}
$$
Surjectivity of $\eee^*$ now follows from that of $\alpha_x^*$ and from proposition \ref{alphabeta}.
\qed

\begin{cor}
\label{Qdxdualagain}
If $d=1$ or 2, then
for each $x\in J^d$ we have (compare with (\ref{Qdxdual}) in section \ref{abeljacobi}
above):
$$
\bp Q_d(x)\dual = {\rm Span}(\kum(W_d -x)) = {\rm Span}(\phi \circ\eeetil \ \widetilde{\bp}(x)) \subset |2\Theta|.
$$
\end{cor}

\pf
We have $\bp Q_d(x)\dual = {\rm Span}(\kum(W_d -x))$ by construction; the second equality then follows from diagram (\ref{abdiagram}), proposition \ref{alphabeta} and corollary \ref{completeness}. Surjectivity of $\alpha_x^*$, when $d=1$ or 2, follows from proposition \ref{conjOK}. 
\qed

In the case $d=1$ we have $\widetilde\bp (x) = \bp (x)$, and this extension space embeds linearly under $\phi \circ \eee$. It follows from this that:

\begin{cor}
\label{G2}
$G_2 \subset \phi(\su(2))$.
\end{cor}

\section{Prym quadrisecant planes are cut out by $G_3$}
\label{prymG3}

In \cite{OPP} it was shown that the ruled subvariety $G_2 \subset \phi(\su(2))$ contains the Fay trisecant lines of the Kummer variety in the planes of its ruling. In this section we shall show that in a curiously analogous way the ruling of $G_3 \subset |2\Theta|$ cuts out the Beauville--Debarre quadrisecant planes of each Prym Kummer variety (see section \ref{quadrisecants}).

For $x\in J^2_C$ let us write 
$$
W(x) = \bp Q_2 (x)\dual = \{\ D\in |2\Theta|\ | \ x+W_{g-3} \subset {\rm supp}\ D\ \ \};
$$
equivalent descriptions being given by corollary \ref{Qdxdualagain}. Then set-theoretically 
$
G_3 = \bigcup_{x\in J_C^2} W(x)
$
(see (\ref{NQsequence}) and (\ref{Qdxdual})).
For $\eta\in J_C[2]\backslash\{0\}$ 
we fix $\zeta \in J_C$ such that $\zeta^2 = \eta$; with respect to this $\zeta$ the Prym-Kummer map is then described in diagram (\ref{prymkum}) in section \ref{pryms}. The description of the quadrisecant planes of the image that we shall use is that given by the map $f:\ff \rightarrow S^4 \kum(\peta)$ (see diagram (\ref{FBDdiagram})) where $\ff = \{(\lambda, D)\in J^2 \times S^4 C\,|\, D\in |\lambda^2|\}$. 

\begin{prop}
\label{WcapPrym}
For each $x\in J_C^2$ we have
$$
W(x) \cap \kum(\peta) = \bigcup_{D\in |x^2 \eta|} f(x\zeta, D).
$$
\end{prop}

\pf
By construction, if $u\in \peta$ then $\kum(u) \in W(x)$ if and only if
$ h^0(C,x(D) \otimes \zeta \otimes \pi_* u) >0$ for all $D\in S^{g-3} C$. 
By (\ref{killingD}) this is equivalent to 
$$
0< h^0(C,x \otimes \zeta \otimes \pi_* u) = h^0(\ctil, u\otimes \pi^*(x \zeta)),
$$
and since $\deg u\otimes \pi_*(x \zeta) = 4$ we conclude that $\kum(u) \in W(x)$ if and only if $u = \pi^*(x^{-1} \zeta^{-1})\otimes \oo(\Gamma)$ for some $\Gamma \in S^4 \ctil$. 

By definition of $\peta \cup \peta^-$ we have $\nm_{\pi}(u) = \oo_C$ and so 
$\nm_{\pi} (\Gamma) = D$ for some $D\in |(x \zeta)^2 | = |x^2 \eta|$. Conversely, for each such effective divisor $D$ we have sixteen possibilities for $\Gamma$, giving the sixteen points of $\peta \cup \peta^-$ in table (\ref{FBDtable}) (with $\lambda = x\zeta$).
\qed

Finally, let us restrict attention to 
a nonhyperelliptic curve $C$ of genus 4. We arrive at a configuration
$$
G_2 \subset \phi(\su(2)) \subset G_3 \subset |2\Theta| = \bp^{15}.
$$
Here $G_2$ is ruled by $\bp N_2 \cong \bp Q_1\dual \rightarrow J^1$ with 4-dimensional projective fibres which we shall denote by 
$$
\bp^4(y) = \{D\in |2\Theta|\,|\,y+W_2 \subset {\rm supp}\,D\}
\qquad
\hbox{for $y\in J^1$;}
$$
and $G_3$ is ruled by  
$\bp N_1 \cong \bp Q_2\dual  \rightarrow J^2$ with 10-dimensional fibres
$$
 W(x) = \bp^{10}(x) = \{D\in |2\Theta|\,|\,x+W_1 \subset {\rm supp}\,D\}
\qquad
\hbox{for $x\in J^2$.}
$$
With respect to any theta characteristic $\kappa \in \vartheta(C)$ the two rulings are polar; equivalently they 
are given by the same grassmannian map (where both vertical maps are isomorphisms):
$$
\begin{array}{ccc}
J^1 & \map{}& \gr_{5} {H^0(2\Theta)}\\
&&\\
\downarrow \kappa &&\downarrow {\rm ann}\\
&&\\
J^2 &\map{}& \gr_{11} {H^0(2\Theta)}.\\
\end{array}
$$

We wish to consider the restriction of these rulings to the fixed-point set $|2\Theta|^{\eta}  \subset |2\Theta|$ under the action of a 2-torsion point  
$\eta\in J[2]\backslash \{0\}$.
We recall from section \ref{pryms} 
that this set has two components
$$
|2\Theta|^{\eta} = \bp (V\dual)^+_{\eta} \cup \bp (V\dual)^-_{\eta}
$$
into which $\peta \cup \peta^-$ map (given a choice of $\zeta = {1\over 2}\eta$)
as Kummer varieties. 
We shall identify $\bp (V\dual)^{\pm}_{\eta} = |2\Xi|^{\pm}$ and view $|2\Xi|^+ \cup |2\Xi|^- \subset |2\Theta|$.

\begin{lemm}
\label{intplane}
For $x\in J^2$ we have 
$$
\dim \bp^{10}(x) \cap |2\Xi|^{\pm} = 
\cases{2& if $h^0(x^{2}\eta) =1$\cr
       3& if $h^0(x^{2}\eta) =2$.\cr}
$$
\end{lemm}

\pf
By corollary \ref{intdims} and Serre duality this is equivalent, for any choice of $\kappa \in \vartheta(C)$, to 
$$
\bp^{4}(y) \cap |2\Xi|^{\pm} = 
\cases{\emptyset& if $h^0(y^{2}\eta) =0$\cr
       {\rm point}& if $h^0(y^{2}\eta) >0$\cr}
$$
where $y = \kappa x^{-1}$. Here $\bp^4(y)\hookrightarrow \su(2)$ is the set of extensions $\ses{y^{-1}}{E}{y}$, and $E\in \bp^4(y)$ is in $|2\Theta|^{\eta}$ 
if and only if $E\in \bp^4(y) \cap \bp^4(\eta y)$. 
We now refer to \cite{OPP}, proposition 1.2. This says, if we translate from $\su(2,K)$ to $\su(2)$, that $\bp^4(y) \cap \bp^4(y')$ is nonempty only if either $y\otimes y' = \oo(p+q)$ or $y^{-1}\otimes y' = \oo(p-q)$ for some points $p,q \in C$. When $y' = y \eta$ the second case is impossible for nonhyperelliptic $C$; so we see that the intersection is nonempty only in the first case, i.e. if and only if $h^0(y^{2}\eta) >0$. 
In this case $\bp^4(y) \cap \bp^4(y \eta)$ is a line on which the involution $\eta$ acts with two fixed points, which are the 
respective intersections 
$\bp^{4}(y) \cap |2\Xi|^{\pm}$.
\qed

We conclude from  \ref{WcapPrym} and \ref{intplane}
that the 10-planes $\bp^{10}(x) \subset |2\Theta|$ cut out in $|2\Xi|$ {\it precisely the quadrisecant planes of the embedded Prym-Kummer variety}---and this concludes the proof of theorem \ref{theorem2}.

\section{Appendix: symmetric products of a curve}
\label{symm}

We shall gather together here some results concerning the symmetric products $S^d C$. We denote by 
$
p: C^d = C\times \cdots \times C \rightarrow S^d C
$
the quotient map. For any $L \in \pic(C)$ one can associate a line bundle $S^d L \in \pic\ S^d C$ satisfying 
$$
p^* S^d L = \bigoplus_{i=1}^d pr_i^* L.
$$
Let $\Delta $ be the union of 
the diagonals in $C^d$, i.e. the ramification divisor of $p$; and $\delbar$ the diagonal divisor in $S^d C$; so we have 

\begin{equation}
\label{pbdiag}
2\Delta = p^* \delbar.
\end{equation}
Noting that $K_{C^d} = p^* S^d K_C$, the ramification formula for the map $p$ tells us that 
$$
\oo(\Delta) = p^*(S^d K_C \otimes K^{-1}_{S^d C}).
$$
In other words, although the divisor $\Delta$ does not descend to the quotient, the line bundle $\oo(\Delta)$ does descend. 

We next introduce some notation for the cohomology ring of $S^d C$ (see \cite{Mac}).
Let $\beta \in H^2(C,\z)$ be the fundamental class of a point; and let $\{\alpha_i\}_{1\leq i \leq 2g}$ be a symplectic basis of $H^1(C, \z)$. Then it is well-known that the cohomology ring of $S^d C$ is generated over $\z$ by the following elements (see \cite{Mac} (3.1) and (6.3); and note that we identify $H^*(S^d C, \z)$ with invariant cohomology on $C^d$ under the action of the symmetric group):
$$
\begin{array}{rcll}
\xi_i &=& \displaystyle
\sum_{j=1}^d 1\otimes \cdots 
\underbrace{\otimes \alpha_i \otimes}_{\hbox{\sixrm j-th positions}} 
\cdots \otimes 1 &\in H^1(S^d C, \z) 
\qquad
\hbox{for 
$i = 1,\ldots , 2g$;}\\
\eta &=&\displaystyle \sum_{j=1}^d 1\otimes \cdots 
\overbrace{\otimes \beta \otimes}
\cdots \otimes 1 &\in H^2(S^d C,\z).\\
\end{array}
$$
The $\xi_i$s anticommute with each other and commute with $\eta$. In fact $\eta$ is the fundamental class of the divisor $S^{d-1} C \hookrightarrow S^d C$, $D\mapsto D+p$, where $p\in C$ is any point of the curve. More generally one has

\begin{equation}
\label{c1SdL}
c_1(S^d L) = (\deg L) \eta
\qquad
\hbox{for any $L\in \pic(C)$;}
\end{equation}
while for the inclusion $S^{d-1} C \hookrightarrow S^d C$ above, 
$\oo_{S^d C}(S^{d-1} C) = S^d \oo_C(p)$. We shall sometimes (notably in proposition
\ref{100} below) simply write $p$ for this line bundle.

Next, we define
$$
\sigma_i = \xi_i \xi_{i+g} \in H^2(S^d C, \z)
\qquad
\hbox{for 
$i = 1,\ldots , g$.}
$$
Then we have (\cite{Mac} 5.4): 

\begin{equation}
\label{50}
\begin{array}{rcll}
\sigma_i ^2 &=& 0,&\\
\sigma_i \sigma_j &=& \sigma_j \sigma_i &\hbox{for $i\not= j$,}\\
\sigma_{i_1} \cdots \sigma_{i_a} \eta^b &=& \eta^{a+b}&
\hbox{for $b>0$ and distinct ${i_1}, \ldots ,{i_a}$;}\\
\end{array}
\end{equation}
and the Chern class of the diagonal is:

\begin{equation}
\label{c1delbar}
c_1(\delbar) = 2(d+g-1)\eta -2(\sigma_1 + \cdots + \sigma_g).
\end{equation}

We now consider, for $x\in J^d$, the map

\begin{equation}
\label{alphax}
\begin{array}{rcl}
\alpha_x : S^d C &\rightarrow& J\\
D &\mapsto& x^{-1} \otimes \oo(D).\\
\end{array}
\end{equation}
We shall write $\ll_x = \alpha_x ^* \ll \in \pic\ S^d C$.

\begin{prop}
\label{100}
\begin{enumerate}
\item
$\ll_x \cong S^d (K x^2) \otimes \oo(-\delbar)$;
\item
$c_1(\ll_x) = 2(\sigma_1 + \cdots + \sigma_g)$;
\item
$\chi (S^d C, \ll_x) = \sum_{i=0}^d {g \choose i}$;
\item
$\chi (S^d C, \ll_x(-p)) = {g \choose d}$ for any $p\in C$.
\end{enumerate}
\end{prop}

\pf
1. 
The case $d=1$ is well-known; see for example \cite{NR2}. For $d>1$ let $y=x(-p_1-\cdots -p_{d-1})$ for some arbitrary points $p_1,\ldots ,p_{d-1}\in C$. Then, from the case $d=1$, we see that 
$$
p^*\ll_x|_{\{(p_1,\ldots ,p_{d-1})\}\times C} = Kx^2 (-2p_1-\cdots -2p_{d-1}),
$$
and hence $p^* \ll_x = p^*S^d(Kx^2) \otimes \oo(-2\delbar)$. So the result follows from (\ref{pbdiag}) and injectivity of $p^*$ on line bundles.

2. 
This is immediate from part 1, (\ref{c1SdL}) and (\ref{c1delbar}). 

3. 
One should be able to extract this from \cite{Mac}, but as there appears to be an error in the proof of (14.9) of that paper we shall give the Riemann-Roch calculation here.

By \cite{Mac}(14.5) the total Chern class is
$
c(S^d C) = (1+\eta)^{d-2g+1} \prod_{i=1}^g (1+\eta -\sigma_i),
$
and the Todd class is therefore
$$
\td (S^d C) = \Bigl({\eta \over 1-e^{-\eta}}\Bigr)^{d-2g+1}
             \prod_{i=1}^g 
              \Bigl({\eta -\sigma_i \over 1-e^{-\eta +\sigma_i }}\Bigr).
$$
Making use of (\ref{50}), the expression in the product can be rewritten:
$$
{\eta -\sigma_i \over 1-e^{-\eta +\sigma_i }}
= 
\Bigl({\eta \over 1-e^{-\eta}}\Bigr)(1+\sigma_i \tau)
\qquad
{\rm where}
\ \tau = {\eta e^{-\eta} + e^{-\eta} - 1
        \over \eta (1-e^{-\eta})},
$$
and so
$$
\td (S^d C) = \Bigl({\eta \over 1-e^{-\eta}}\Bigr)^{d-g+1}
                \prod_{i=1}^g (1+\sigma_i \tau).
$$
On the other hand ${\rm ch}(\ll_x) = e^{2(\sigma_1 + \cdots + \sigma_g)}
= \prod_{i=1}^g (1+2 \sigma_i)$; 
so by Hirzebruch-Riemann-Roch:
$$
\chi(\ll_x) = \deg \bigg\lbrace
\Bigl({\eta \over 1-e^{-\eta}}\Bigr)^{d-g+1}
                \prod_{i=1}^g (1+\sigma_i (2+\tau))
\bigg\rbrace _d
$$
Using \cite{Mac} (5.4) this expression may be computed by setting $\sigma_1 =\cdots = \sigma_g = \eta$ and evaluating the coefficient of $\eta^d$:
$$
\begin{array}{rcl}
\chi(\ll_x) &=& 
 \bigg\lbrace                     
\Bigl({\displaystyle\eta \over \displaystyle 1-e^{-\eta}}\Bigr)^{d-g+1}(1+\eta (2+\tau))^g
\bigg\rbrace _{\hbox{coefficient of $\eta^d$}}
\\
&=&
 \bigg\lbrace
\Bigl({\displaystyle\eta \over \displaystyle 1-e^{-\eta}}\Bigr)^{d+1}(2-e^{-\eta})^g
\bigg\rbrace _{\hbox{coefficient of $\eta^d$}}\\
&&\\
&=&\displaystyle
{\rm Res}_{\eta =0}\ {(2-e^{-\eta})^g d\eta\over (1-e^{-\eta})^{d+1}} \\
&&\\
&=&\displaystyle
{\rm Res}_{\zeta =0}\ {(1+\zeta)^g d\zeta\over \zeta^{d+1}(1-\zeta)}\\
&=&\displaystyle
\sum_{i=0}^d {g \choose i},\\
\end{array}
$$
where in the penultimate step we have made a substitution $\zeta = 1 - e^{-\eta}$.

The computation for part 4 is entirely similar.
\qed

For $x\in J^d$ ($d\geq 1$), let $\bp (x) \cong \bp ^{g+2d-2}$ 
be the extension space introduced in section \ref{segrestrat}
into which $C$ maps by the complete linear series $|Kx^2|$:

\begin{equation}
\lambda_{|Kx^2|}: C \hookrightarrow \bp (x) = |Kx^2|\dual.
\end{equation}
We shall denote by $\iic$ the ideal sheaf of the image, and consider the linear series $|\iic^{d-1}(d)|$ on $\bp(x)$. This series contracts each $d$-secant plane $\dbar \cong \bp^{d-1} \subset \bp(x)$ (where $D\in S^d C$) to a point, and so induces a rational map 

\begin{equation}
\label{betax}
\beta_x : S^d C \rightarrow |\iic^{d-1}(d)|\dual.
\end{equation}

We observe that by proposition \ref{100}(1) a section of $\ll_x$ can be identified with a symmetric multilinear form 
$F: \bigotimes^d H^0(C, Kx^2)\dual \rightarrow \c$
with the property that $F(p_1,\ldots , p_d) = 0$ whenever $p_i = p_j$, $i\not= j$; or equivalently $F \in H^0(\bp(x), \iic^{d-1}(d))$. From this one may check the following (and we are grateful to Beauville for pointing this out to us):

\begin{prop}
\label{alphabeta}
$\beta_x^* \oo(1) = \ll_x$; and there is an induced isomorphism
$$
\beta^*_x : H^0(\bp (x), \iic^{d-1}(d)) \idfy  H^0(S^d C, \ll_x).  
$$
\end{prop}

One would like higher cohomology of $\ll_x$ to vanish, so that parts 3 and 4 of proposition \ref{100}
compute $h^0$.

\begin{conj}
\label{400}
Let $1\leq d\leq g$ and $x\in J_C^d$. Then:
\begin{enumerate}
\item
$h^0(S^d C, \ll_x) = \sum_{i=0}^d {g\choose i}$;
\item
$\alpha_x^* : H^0(J,\ll) \rightarrow H^0(S^d C, \ll_x)$ is surjective.
\end{enumerate}
\end{conj}

We can verify this at once in certain cases. When $d=1$, for example, $\ll_x = K x^2$ and \ref{400} follows from \cite{NR2} lemma 4.1. In the case $d=2$, part 1 is proved in \cite{BV} proposition 4.9. Part~2 can also be shown---for $g>4$ this will appear elsewhere, but in the case $g= 4$ we can give a simple ad hoc argument as follows. In the notation of section \ref{abeljacobi}, $N_2(x) \perp N_1(\kappa x^{-1})$ where $\dim N_1(\kappa x^{-1}) = 11$. Hence $\dim N_2(x) \leq 5$; but (by \cite{BV} proposition 4.9) $\dim Q_2(x) = 11$, and therefore $\alpha_x^*$ is surjective.

In the case $d=g$ the conjecture is also easy to see, since $\alpha_x : S^g C \rightarrow J_C$ is a birational morphism with connected fibres. Indeed, this suggests using a descending induction for $d\leq g$ using the inclusions 
$$
C\subset S^2 C \subset \cdots \subset S^g C
$$
given by some choice of point $p\in C$. At each stage restriction yields an exact sequence on $S^{d}C$ (where the last term is supported on $S^{d-1} C$):
$$
\ses{\ll_{x}(-p)}{\ll_{x}}{\ll_{x(-p)}}.
$$

Suppose for a moment that 

\begin{equation}
\label{vanishing}
H^i(S^{d} C , \ll_{x} (-p))=0
\qquad
\hbox{for $i>0$.}
\end{equation} 
This
would imply that the restriction map 
$
H^0(S^{d} C , \ll_{x} ) \rightarrow H^0(S^{d-1} C , \ll_{x(-p)} )
$
is surjective, and hence surjectivity of $\alpha_{x}^*$ would imply that of $\alpha_{x(-p)}^*$; while \ref{400}(1) also follows inductively since $h^0(S^{d} C , \ll_{x} (-p)) = {g\choose d}$ by proposition \ref{100}(4).

Unfortunately we cannot prove (\ref{vanishing})---except when $d=g$. For here we have a diagram

$$
\begin{array}{rcl}
S^{g-1} C &\subset&S^{g} C\\
\downarrow &&\downarrow \\
\Theta &\subset& J\\
\end{array}
$$
where the vertical arrows are birational morphisms with connected fibres and so preserve cohomology under pull-back. Hence $h^i(S^g C,\ll_{x} (-p)) = h^i(J,\ll(-\Theta)) = h^i(J,\Theta) = 0$.

To summarise, we have shown:

\begin{prop}
\label{conjOK}
Conjecture \ref{400} holds for $d=1,2,g-1,g$; and in particular it holds when~$g\leq 4$.
\end{prop}

Finally, we remark that the difficulty in proving (\ref{vanishing}) rests 
principally in the fact that the self-intersection
$$
c_1(\ll_x K^{-1})^d = {g!\over (g-d)!} + (d+1)^d - g^d
$$
is negative in all the cases not covered by proposition \ref{conjOK}, 
and one can show that the line bundle $\ll_x K^{-1}$ fails to be nef, big or semi-ample. One is therefore unable to use the standard vanishing theorems.

\bigskip

\noindent
{\addressit
Department of Mathematical Sciences, Science Laboratories, South Road, Durham DH1~3LE, U.K.
}

\noindent
{\addressit E-mail:} {\eightrm w.m.oxbury@durham.ac.uk}

\medskip
\noindent
{\addressit
DPMMS, University of Cambridge, 16 Mill lane, Cambridge CB2~1SB, U.K.
}

\noindent
{\addressit E-mail:} {\eightrm pauly@pmms.cam.ac.uk}

\end{document}